\documentclass[a4paper,11pt,reqno,noindent]{amsart}
\usepackage[centertags]{amsmath}
\usepackage{amsfonts,amssymb,amsthm,dsfont,cases,amscd,esint,enumerate}
\usepackage[T1]{fontenc}
\usepackage[english]{babel}
\usepackage[applemac]{inputenc}
\usepackage{newlfont}
\usepackage{color}
\usepackage[body={15cm,21.5cm},centering]{geometry} 
\usepackage{fancyhdr}
\usepackage{enumerate}

\theoremstyle{plain}
\newtheorem{theor10}{Theorem}

\newtheorem{prop10}{Proposition}

\newtheorem{cor10}{Corollary}

\newtheorem{theor0}{Theorem}[section]

\newtheorem{lem0}[theor0]{Lemma}

\newtheorem{prop0}[theor0]{Proposition}

\newtheorem{cor0}[theor0]{Corollary}

\newtheorem{propr0}[theor0]{Property}

\newtheorem{hyp0}[theor0]{Hypothesis}

\newtheorem{result0}[theor0]{Result}

\newtheorem{conj0}[theor0]{Conjecture}

\newtheorem{heur0}[theor0]{Heuristics}

\theoremstyle{definition}
\newtheorem{defin0}[theor0]{defin0}

\newtheorem{rems0}[theor0]{Remark}

\newtheorem{ex0}[theor0]{Example}

\newtheorem{exs0}[theor0]{Examples}

\newtheorem{rem0}[theor0]{Remark}

\newtheorem{qu0}[theor0]{Question}

\newtheorem{qus0}[theor0]{Questions}

  \newtheorem{as0}[theor0]{Assumption}

\mathchardef\emptyset="001F
\numberwithin{equation}{section}

\usepackage{color}
\usepackage[colorlinks,citecolor=black,urlcolor=black]{hyperref}

\title[An algebraic approach for solving fourth-order PDFs]{An algebraic approach for solving fourth-order partial differential equations}

\author{
A.~ Pogorui, 
T.~ Kolomiiets 
and R.~M.~ Rodr\'{\i}guez-Dagnino 
}

\begin{document}

\begin{abstract}
It is well-known that any solution of the Laplace equation is a real or imaginary part of a complex holomorphic function.
In this paper, in some sense, we extend this property into four order hyperbolic and elliptic type PDEs.
To be more specific, the extension is for a $c$-biwave PDE with constant coefficients, and we show that the components of a differentiable function on the associated hypercomplex algebras provide solutions for the equation.
\end{abstract}



\maketitle

\section{Introduction}

In this paper we are interested in finding the solution of the following equation
\begin{equation}
\left(\frac{{\partial }^4}{\partial x^4}-2 c \frac{{\partial }^4}{\partial x^2\partial y^2}+\frac{{\partial }^4}{\partial y^4}\right)u\left(x,y\right) = 0,\, \, c>0.  \label{GrindEQ__1_}
\end{equation}
Depending on the value of $c$ we may consider three cases. Namely, the case where $0<c<1$ and we call it as the
$c$-biwave equation of the elliptic type, the case where $c>1$ and we call it as the $c$-biwave equation of hyperbolic type,
and in the case where $c=1$ Eq.\eqref{GrindEQ__1_} is the well-known biwave equation.
The biwave equation has been used in modeling of
$d$-wave superconductors (see for instance \cite{Feng11}, and references therein) or in probability
theory \cite{Pogor05,Kolo15}. In  \cite{Grysh18} the author studied Eq.\eqref{GrindEQ__1_} in the case where $c<-1$ and
considered its application to theory of plain orthotropy.

It is easily verified that any equation of the form
\[\left(A\frac{{\partial }^4}{\partial x^4}+2B\frac{{\partial }^4}{\partial x^2\partial y^2}+C\frac{{\partial }^4}{\partial y^4}\right)u\left(x,y\right)=0,\]
where $AC>0$ and $AB<0$ can be reduced to Eq.\eqref{GrindEQ__1_} by changing variables. To obtain all solutions of Eq. \eqref{GrindEQ__1_} for $1\ne c>0$ we will use the method developed in \cite{Pogor14}. According to such approach we need a commutative algebra with basis containing $e_1$, $e_2$ such that
\begin{equation}
e^4_1-2 c\,e^2_1 e^2_2+e^4_2=0. \label{GrindEQ__2_}
\end{equation}

Then, we study monogenic functions on the subspace of this algebra containing $e_1$, $e_2$ and show that any solution
of Eq. \eqref{GrindEQ__1_} can be obtained as a component of such monogenic functions.

\section{Hyperbolic case}

Firstly we study Eq. \eqref{GrindEQ__1_} in the case where $c>1$, which is said to be hyperbolic. Let us consider an
associative commutative algebra over the real field ${\mathbb R}$
$$A_c=\left\{x\mathbf {u} +y \mathbf {f}+z \mathbf {e}+v \mathbf {fe} : x,y,z,v \in \mathbb{R} \right\}$$
  with a basis ${\mathbf u}$, ${\mathbf f}$,
${\mathbf e}$, $\mathbf {fe}$, where ${\mathbf u}$ is the identity element of $A_c$ and the following Cayley table holds
$\mathbf {fe}=\mathbf {ef}$,
${\mathbf {f}}^2 = {\mathbf {u}}$, ${{\mathbf e}}^2 = {\mathbf u}-m \mathbf {fe}$, where $m=\sqrt{2(c-1)}$.

The basis elements ${\mathbf u}$, ${\mathbf e}$ satisfy Eq. \eqref{GrindEQ__2_}.

It is easily verified that for $c>1$ algebra $A_c$ has the following idempotents
\[
i_1=\frac{k_1}{k_1+k_2}{\mathbf u}-\frac{\mathbf {f}\sqrt{2}}{k_1+k_2}{\mathbf e},
\]
\begin{equation}
i_2=\frac{k_2}{k_1+k_2}{\mathbf u}+\frac{\mathbf {f}\sqrt{2}}{k_1+k_2}{\mathbf e}, \label{GrindEQ__3_}
\end{equation}
where $k_1=\sqrt{c+1}-\sqrt{c-1}$, $k_2=\sqrt{c+1}+\sqrt{c-1}$.

Therefore, we have
\[i_1+i_2= {\mathbf u}\]
and

\begin{align*}
i_1 \,i_2=\frac{k_1 k_2}{{\left(k_1+k_2\right)}^2}{\mathbf u}-\frac{\sqrt{2}k_2}{{\left(k_1+k_2\right)}^2}\mathbf {f}{\mathbf e}
+\frac{\sqrt{2}k_1}{{\left(k_1+k_2\right)}^2}\mathbf {f}{\mathbf e} \\
-\frac{2}{{\left(k_1+k_2\right)}^2}{\mathbf u}
+\frac{2m}{{\left(k_1+k_2\right)}^2}\mathbf {f}{\mathbf e}=0.
\end{align*}

It is easily seen that
\begin{equation}
{\mathbf e} ={\mathbf f} \frac{k_1}{\sqrt{2}}i_2 - {\mathbf f} \frac{k_2}{\sqrt{2}}i_1. \label{GrindEQ__4_}
\end{equation}

Consider a subspace $B_c$ of algebra $A_c$ of the following form
\[B_c=\left\{x{\mathbf u}+y{\mathbf e}\mathrel{\left|\vphantom{x{\mathbf u}+y{\mathbf e} x,y\in {\mathbb R}}\right.\kern-\nulldelimiterspace}x,y\in {\mathbb R}\right\}.\]

\begin{defin0} \textit{A function } $g:\,B_c\to A_c$ \textit{ is called differentiable (or monogenic) on }
$B_c$ \textit{ if for any } $B_c\ni w=x{\mathbf u}+y{\mathbf e}$ \textit{ there exists a unique element }
$g'\left(w\right)$ \textit{ such that for any } $h\in B_c$

\[\lim_{{\mathbb R}\ni \varepsilon \to 0} \frac{g\left(w+\varepsilon h\right)-g\left(w\right)}{\varepsilon }
=h g'\left(w\right),\]
\textit{where }$h g'\left(w\right)$\textit{ is the product of }$h$\textit{ and }$g'\left(w\right)$\textit{ as elements of} $A_c$.
\end{defin0}

It follows from \cite{Pogor14} that a function $g\left(w\right)={\mathbf u}\,u_1\left(x,y\right)+{\mathbf f}\,u_2\left(x,y\right)
+{\mathbf e}\,u_3\left(x,y\right)+{\mathbf {f\,e}} \,u_4\left(x,y\right)$
is monogenic if and only if there exist continuous partial derivatives $\frac{\partial u_i\left(x,y\right)}{\partial x}$, $\frac{\partial u_i\left(x,y\right)}{\partial y},$ $i=1,2,3,4$ and it satisfies the following Cauchy-Riemann type conditions
\begin{equation*}
\, {\mathbf e}\frac{\partial }{\partial x}g\left(w\right)={\mathbf u}\frac{\partial }{\partial y}g\left(w\right), \, \forall w\in B_c,
\end{equation*}
or
\[\frac{\partial u_1\left(x,y\right)}{\partial y}=\frac{\partial u_3\left(x,y\right)}{\partial x},\]
\[\frac{\partial u_2\left(x,y\right)}{\partial y}=\frac{\partial u_4\left(x,y\right)}{\partial x},\]
\[\frac{\partial u_3\left(x,y\right)}{\partial y}=\frac{\partial u_1\left(x,y\right)}{\partial x}-m\frac{\partial u_4\left(x,y\right)}{\partial x},\]
\[\frac{\partial u_4\left(x,y\right)}{\partial y}=\frac{\partial u_2\left(x,y\right)}{\partial x}-m\frac{\partial u_3\left(x,y\right)}{\partial x}.\]

It is also proved in \cite{Pogor14} that if $g$
is monogenic then its components $u_i\left(x,y\right)$ satisfies Eq. \eqref{GrindEQ__1_}.

By passing in $B_c$ from the basis ${\mathbf u}$, ${\mathbf e}$ to the basis $i_1$, $i_2$, we have
\[w=x\,{\mathbf u}+y\,{\mathbf e}=\left(x-{\mathbf f}\frac{k_2}{\sqrt{2}}y\right)\,i_1 +\left(x+{\mathbf f}\frac{k_1}{\sqrt{2}}y\right)\,i_2.\]

\begin{lem0}
\textit{A function }$g:\,B_c\to A_c$\textit{, where }$c>1$\textit{, is differentiable if and only if it can be
represented as follows}
\begin{equation} g\left(w\right)=\alpha \left(w_1\right)i_1+\beta \left(w_2\right) i_2, \label{GrindEQ__15_}
\end{equation}
\textit{where }$w_1=x-{\mathbf f}\frac{k_2}{\sqrt{2}}y$\textit{, }$w_2=x+{\mathbf f}\frac{k_1}{\sqrt{2}}y$
\textit{ and }$\alpha \left(w_1\right)$\textit{, }$\beta \left(w_2\right)$
\textit{ have continuous partial derivatives $\frac{\partial }{\partial x}\alpha \left(w_1\right), \frac{\partial }{\partial y}\alpha \left(w_1\right), \frac{\partial }{\partial x}\beta \left(w_2\right), \frac{\partial }{\partial y}\beta \left(w_2\right)$ satisfying }
\[\frac{\partial }{\partial y}\alpha \left(w_1\right)=-{\mathbf f}\frac{k_2}{\sqrt{2}} \frac{\partial }{\partial x}\alpha \left(w_1\right),\]  \[\frac{\partial }{\partial y}\beta \left(w_2\right)={\mathbf f}\frac{k_1}{\sqrt{2}} \frac{\partial }{\partial x}\beta \left(w_2\right).\]
\end{lem0}
\begin{proof}
Sufficiency can be verified directly. Indeed,
\[\frac{\partial }{\partial y}g\left(w\right)=\frac{\partial }{\partial y}\alpha \left(w_1\right)i_1+\frac{\partial }{\partial y}\beta \left(w_2\right)i_2\]
\[=-{\mathbf f}\frac{k_2}{\sqrt{2}} \frac{\partial }{\partial x}\alpha \left(w_1\right)i_1+{\mathbf f}\frac{k_1}{\sqrt{2}} \frac{\partial }{\partial x}\beta \left(w_2\right)i_2\]

On the other hand, taking into account Eqs. \eqref{GrindEQ__3_}, \eqref{GrindEQ__4_}, we have
\[{\mathbf e}\frac{\partial }{\partial x}g\left(w\right)=\left({\mathbf f}\frac{k_1}{\sqrt{2}}i_2 - {\mathbf f}\frac{k_2}{\sqrt{2}}i_1\right)\left(\frac{\partial }{\partial x}\alpha \left(w_1\right)i_1 + \frac{\partial }{\partial x}\beta \left(w_2\right)i_2\right)\]
\[= -{\mathbf f}\frac{k_2}{\sqrt{2}} \frac{\partial }{\partial x}\alpha \left(w_1\right)i_1 +{\mathbf f}\frac{k_1}{\sqrt{2}} \frac{\partial }{\partial x}\beta \left(w_2\right)i_2.\]

Hence,
\[{\mathbf e}\frac{\partial }{\partial x}g\left(w\right)={\mathbf u}\frac{\partial }{\partial y}g\left(w\right).\]

Now let us prove necessity. Suppose that a function
\[g\left(w\right)={\mathbf u}u_1\left(x,y\right)+{\mathbf f}\,u_2\left(x,y\right)+{\mathbf e}\,u_3\left(x,y\right)
+{\mathbf {f\,e}}\,u_4\left(x,y\right)\]
is monogenic on $B_c$. Let us define
\[\alpha \left(w_1\right)={\mathbf u}\left(u_1\left(x,y\right) - \frac{k_2}{\sqrt{2}}u_4\left(x,y\right)\right)+{\mathbf f}\left(u_2\left(x,y\right)-\frac{k_2}{\sqrt{2}}u_3\left(x,y\right)\right),\]
\[\beta \left(w_2\right)={\mathbf u}\left(u_1\left(x,y\right) + \frac{k_1}{\sqrt{2}}u_4\left(x,y\right)\right)+{\mathbf f}\left(u_2\left(x,y\right)+\frac{k_1}{\sqrt{2}}u_3\left(x,y\right)\right).\]

Thus, we have
\[\frac{\partial }{\partial y}\alpha \left(w_1\right)={\mathbf u}\left(\frac{\partial u_3\left(x,y\right)}{\partial x}-\frac{k_2}{\sqrt{2}}\left(\frac{\partial u_2\left(x,y\right)}{\partial x}-m\frac{\partial u_3\left(x,y\right)}{\partial x}\right)\right)\]
\[\qquad + {\mathbf f}\left(\frac{\partial u_4\left(x,y\right)}{\partial x}-\frac{k_2}{\sqrt{2}}\left(\frac{\partial u_1\left(x,y\right)}{\partial x}-m\frac{\partial u_4\left(x,y\right)}{\partial x}\right)\right)\]

\[=-{\mathbf f}\frac{k_2}{\sqrt{2}}\frac{\partial u_1\left(x,y\right)}{\partial x} - {\mathbf u}\frac{k_2}{\sqrt{2}}\frac{\partial u_2\left(x,y\right)}{\partial x}+{\mathbf u}\left(\frac{k_2}{\sqrt{2}}\ m+1\right)\frac{\partial u_3\left(x,y\right)}{\partial x}\]
\[\qquad +{\mathbf f}\left(\frac{k_2}{\sqrt{2}}\ m+1\right)\frac{\partial u_4\left(x,y\right)}{\partial x}.\]

Taking into account that
\[\frac{k_2}{\sqrt{2}} m+1=\sqrt{c^2-1}+c=\frac{k^2_2}{2},\]
we have $\frac{\partial }{\partial y}\alpha \left(w_1\right)=-{\mathbf f}\frac{k_2}{\sqrt{2}} \frac{\partial }{\partial x}\alpha \left(w_1\right)$.

Much in the same manner, it can be shown that $\frac{\partial }{\partial y}\beta \left(w_2\right)={\mathbf f}\frac{k_1}{\sqrt{2}} \frac{\partial }{\partial x}\beta \left(w_2\right)$\textit{.}
\end{proof}
\begin{rems0} \textit{Considering variables }$x$\textit{, }$y_1=-\frac{k_2}{\sqrt{2}}y$\textit{ and }$x$\textit{, }$y_2=\frac{k_1}{\sqrt{2}}y,$\textit{ we have }

\begin{equation} \label{GrindEQ__5_}
\frac{\partial }{\partial y_1}\alpha  = {\mathbf f} \frac{\partial }{\partial x}\alpha,
\end{equation}
\[\frac{\partial }{\partial y_2}\beta  = {\mathbf f} \frac{\partial }{\partial x}\beta.\]

\textit{Hence,} \textit{it is easily verified that if the components $\alpha_1$, $\alpha_2$ of} $\alpha \left(w_1\right)=\alpha_1 \left(w_1\right)+{\mathbf f}\alpha_2 \left(w_1\right)$ \textit{have continuous partial derivatives $\frac{{\partial }^2}{\partial x^2}\alpha_k \left(w_1\right)$ and $\frac{{\partial }^2}{\partial y^2_1}\alpha_k \left(w_1\right)$, $k=1,2$ then they satisfy the wave equation}
\[\left(\frac{{\partial }^2}{\partial x^2}-\frac{{\partial }^2}{\partial y^2_1}\right)u\left(x,y_1\right)=0.\]
\textit{Similarly, the components $\beta_1$, $\beta_2$ of} $\beta \left(w_2\right)=\beta_1 \left(w_2\right)+{\mathbf f}\beta_2 \left(w_2\right)$ \textit{satisfy the wave equation}
\[\left(\frac{{\partial }^2}{\partial x^2}-\frac{{\partial }^2}{\partial y^2_2}\right)u\left(x,y_2\right)=0.\]

\end{rems0}

\begin{theor10} \label{the1}
$u\left(x,y\right)$\textit{ is a solution of Eq. \eqref{GrindEQ__1_} for }$c>1$\textit{ if and only if for some }$i,j\in \left\{1,2\right\}$\textit{ it can be represented as follows}

\[u\left(x,y\right)={\alpha }_i\left({\omega }_1\right)+{\beta }_j\left({\omega }_2\right),\]

\textit{where }${\alpha }_i\left({\omega }_1\right),{\beta }_j\left({\omega }_2\right)$\textit{ are four times continuous differentiable components of }$\ \alpha \left({\omega }_1\right)$\textit{ and }$\beta \left({\omega }_2\right)$\textit{ of monogenic function }$g\left(\omega \right)$\textit{ in the decomposition \eqref{GrindEQ__15_} i.e.,}

\[g\left(\omega \right)=\alpha \left({\omega }_1\right)i_1+\beta \left({\omega }_2\right)i_2,\]

\textit{where }$\alpha \left({\omega }_1\right)={\alpha }_1\left({\omega }_1\right)+{\mathbf f}{\alpha }_2\left({\omega }_1\right)$\textit{, }$\beta \left({\omega }_2\right)={\beta }_1\left({\omega }_2\right)+{\mathbf f}{\beta }_2\left({\omega }_2\right)$\textit{ satisfy Eq.} 
\eqref{GrindEQ__5_}.

\end{theor10}
\begin{proof}
 As it was mentioned above $u\left(x,y\right)={\alpha }_i\left({\omega }_1\right)+{\beta }_j\left({\omega }_2\right)$ is a solution for $c>1$ of Eq. \eqref{GrindEQ__1_}.

Now suppose that $u\left(x,y\right)$ is a solution of Eq. \eqref{GrindEQ__1_}. It is easily verified that

\begin{equation}
\left(\frac{{\partial }^4}{\partial x^4}-2c\frac{{\partial }^4}{\partial x^2\partial y^2}+\frac{{\partial }^4}{\partial y^4}\right)u\left(x,y\right)=\left(\frac{{\partial }^2}{\partial x^2}-\frac{{\partial }^2}{\partial y^2_1}\right)\left(\frac{{\partial }^2}{\partial x^2}-\frac{{\partial }^2}{\partial y^2_2}\right)u\left(x,y\right)=0. \label{GrindEQ__11_}
\end{equation}

It is easily seen that Eq. \eqref{GrindEQ__11_} is equivalent to the set of the following systems

\[\left\{ \begin{array}{c}
\left(\frac{{\partial }^2}{\partial x^2}-\frac{{\partial }^2}{\partial y^2_2}\right)u\left(x,y\right)=v_1\left(x,y\right), \\
\left(\frac{{\partial }^2}{\partial x^2}-\frac{{\partial }^2}{\partial y^2_1}\right)v_1\left(x,y\right)=0 \end{array}
\right.\]

or

\[\left\{ \begin{array}{c}
\left(\frac{{\partial }^2}{\partial x^2}-\frac{{\partial }^2}{\partial y^2_1}\right)u\left(x,y\right)=v_2\left(x,y\right), \\
\left(\frac{{\partial }^2}{\partial x^2}-\frac{{\partial }^2}{\partial y^2_2}\right)v_2\left(x,y\right)=0. \end{array}
\right.\]

Let us consider the first system. Since any solution of $\left(\frac{{\partial }^2}{\partial x^2}-\frac{{\partial }^2}{\partial y^2_k}\right)v_k\left(x,y\right)=0$ is of the form $v_k\left(x,y\right)=f_1\left(x+y_k\right)+f_2\left(x-y_k\right)$, where $f_i$, $i=1,2$ are arbitrary twice differentiable functions it follows from the second equation of the system that

\[v_1\left(x,y\right)=f_1\left(x+y_1\right)+f_2\left(x-y_1\right).\]

Thus, the first equation of the system is

\begin{equation}
\left(\frac{{\partial }^2}{\partial x^2}-\frac{{\partial }^2}{\partial y^2_2}\right)u\left(x,y\right)=f_1\left(x+y_1\right)+f_2\left(x-y_1\right). \label{GrindEQ__12_}
\end{equation}

It is easily seen that a partial solution of Eq. \eqref{GrindEQ__12_} is

\begin{align*}U\left(x,y\right)=\frac{k^2_1}{k^2_1-k^2_2}\left(F_1\left(x-\frac{k_2}{\sqrt{2}}y\right)+F_2\left(x+\frac{k_2}{\sqrt{2}}y\right)\right) \\=\frac{k^2_1}{k^2_1-k^2_2}\left(F_1\left(x+y_1\right)+F_2\left(x-y_1\right)\right),\end{align*}

where $F^{''}_k=f_k$, $k=1,2$.

Thus, the general solution of the system is as follows

\[u\left(x,y\right)=g_1\left(x+y_2\right)+g_2\left(x-y_2\right)+\frac{k^2_1}{k^2_1-k^2_2}\left(F_1\left(x+y_1\right)+F_2\left(x-y_1\right)\right)\]

Let us put ${\alpha }_1\left({\omega }_1\right)=\frac{k^2_1}{k^2_1-k^2_2}\left(F_1\left(x+y_1\right)+F_2\left(x-y_1\right)\right)$ and ${\beta }_2\left({\omega }_2\right)=g_1\left(x+y_2\right)+g_2\left(x-y_2\right)$.

Taking into account that $\left(\frac{{\partial }^2}{\partial x^2}-\frac{{\partial }^2}{\partial y^2_1}\right) \left[\frac{k^2_1}{k^2_1-k^2_2}\left(F_1\left(x+y_1\right)+F_2\left(x-y_1\right)\right)\right]=0$ and \\ $\left(\frac{{\partial }^2}{\partial x^2}-\frac{{\partial }^2}{\partial y^2_2}\right)\left[g_1\left(x+y_2\right)+g_2\left(x-y_2\right)\right]=0$ we conclude the proof for the first system.

The case of the second system can be proved similarly.
\end{proof}
\section{Elliptic case}

Now we consider an associative commutative algebra $A_c$, where $0<c<1$, over the complex field ${\mathbb C}$
with a basis ${\mathbf u}$, ${\mathbf e}$ and the following Cayley table
${\mathbf {u\,e}} ={\mathbf {e\,u}} = {\mathbf e}$, ${{\mathbf e}}^2={\mathbf u}+ {\mathbf i}\mu {\mathbf e}$,
where $\mu =\sqrt{2(1-c)}$. The matrix representations of ${\mathbf u}$ and ${\mathbf e}$ are
\[{\mathbf u}=\left( \begin{array}{cc}
1 & 0 \\
0 & 1 \end{array}
\right){\mathbf ,\ \ }{\mathbf e}=\left( \begin{array}{cc}
0 & 1 \\
1 & {\mathbf i}\mu  \end{array}
\right).\]

Hence, we have the following traces of these representations
\[{tr}\left({\mathbf {uu}}\right)=2, \quad {tr}\left({\mathbf {ue}}\right)={\mathbf i}\mu ,\quad  {tr}\left({\mathbf {ee}}\right)=2-{\mu }^2.\]

Since
\[{det}\left( \begin{array}{cc}
{tr}\left({\mathbf {uu}}\right) & {tr}\left({\mathbf {ue}}\right) \\
{tr}\left({\mathbf {ue}}\right) & {tr}\left({\mathbf {ee}}\right) \end{array}
\right)=2(1+c)\ne 0,\]
then, $A_c$ is a semi-simple algebra \cite{Waerden59}.

By following similar steps as in Eq. \eqref{GrindEQ__3_} we can show that for $0<c<1$ algebra $A_c$ has the
following idempotents
\begin{equation}
I_-=\frac{k_1}{k_1+k_2}{\mathbf u}+\frac{\sqrt{2}}{k_1+k_2}{\mathbf e},\quad
I_+=\frac{k_2}{k_1+k_2}{\mathbf u}-\frac{\sqrt{2}}{k_1+k_2}{\mathbf e}, \label{GrindEQ__6_}
\end{equation}
where $k_1=\sqrt{c+1}-{\mathbf i}\sqrt{1-c}$, $k_2=\sqrt{c+1}+{\mathbf i}\sqrt{1-c}$.

It is also easily verified that these idempotents also satisfy
\[I_-+I_+={\mathbf u}\]
and
\[I_-\,I_+=0.\]

It is straightforward to see that
\begin{equation}
{\mathbf e} = \frac{ k_2}{\sqrt{2}}\,I_- - \frac{ k_1}{\sqrt{2}}\,I_+. \label{GrindEQ__7_}
\end{equation}

\begin{lem0}
\textit{All non-zero elements of subspace }$B_c=\left\{x\,{\mathbf u}+y\,{\mathbf e}\mathrel{\left|\vphantom{x\,{\mathbf u}
+y\,{\mathbf e}\,\, x,y\in {\mathbb R}}\right.\kern-\nulldelimiterspace}x,y\in {\mathbb R}\right\}$
\textit{ of algebra }$A_c$\textit{ are invertible, that is, if }$0\ne w\in B_c$\textit{ then there exists} $w^{-1}\in B_c$.
\end{lem0}
\begin{proof} Suppose $w=s\,{\mathbf u}+t\,{\mathbf e}\in B_c$.
Let us show that there exists $w^{-1}=x\,{\mathbf u}+y\,{\mathbf e}$,
$x,y\in {\mathbb R}$ such that $w\,w^{-1}=1$. Indeed, the equation
\[\left(s\,{\mathbf u}+t\,{\mathbf e}\right)\left(x\,{\mathbf u}+y\,{\mathbf e}\right)={\mathbf u}\]
has a unique solution since the determinant of the system
\[ \begin{array}{c}
sx+ty=1, \\
tx+\left(s+{\mathbf i}\mu t\right)y=0, \end{array}
\]
where $x,y$ are unknown, is $\Delta =s^2-t^2+{\mathbf i}\mu ts$ and $\Delta =0$ if and only if $s=t=0$.
A function $f\left(w\right)$, $w\in B_c$ is said to be differentiable if it is differentiable in the common sense, i.e., for all
$w\in B_c$ there exists the following limit
\[\lim_{B_c \ni \Delta w\to 0} \frac{f\left(w+\Delta w\right)-f\left(w\right)}{\Delta w} =f'\left(w\right).\]

It is easily seen that if $f$ is differentiable then it is monogenic and hence, it satisfies the following Cauchy-Riemann type of conditions \cite{Pogor14}
\[
{\mathbf e}\frac{\partial }{\partial x}f\left(w\right)={\mathbf u}\frac{\partial }{\partial y}f\left(w\right)
\]
or in this case we have
\[\frac{\partial u_1\left(x,y\right)}{\partial y}=\frac{\partial u_3\left(x,y\right)}{\partial x},\quad   \frac{\partial u_2\left(x,y\right)}{\partial y}=\frac{\partial u_4\left(x,y\right)}{\partial x},\]

\[\frac{\partial u_3\left(x,y\right)}{\partial y}=\frac{\partial u_1\left(x,y\right)}{\partial x}-\mu \frac{\partial u_4\left(x,y\right)}{\partial x},\]
\[\frac{\partial u_4\left(x,y\right)}{\partial y}=\frac{\partial u_2\left(x,y\right)}{\partial x}+\mu \frac{\partial u_3\left(x,y\right)}{\partial x}.\]

In \cite{Pogor14} it is also proved that if a function $f\left(x,y\right)= {\mathbf u}u_1\left(x,y\right)+{\mathbf i}\,u_2\left(x,y\right)
+{\mathbf e}\,u_3\left(x,y\right)+{\mathbf i}\,{\mathbf e}\,u_4\left(x,y\right)$ is monogenic
then $u_i\left(x,y\right)$ satisfies Eq. \eqref{GrindEQ__1_}. We should mention that a constructive description of
monogenic functions in a three-dimensional harmonic algebra was studied in \cite{Plaksa10,Plaksa13}.

By passing from the basis ${\mathbf u}$, ${\mathbf e}$ to the basis $I_-$, $I_+$ we have
\[w=x\,{\mathbf u}+y\,{\mathbf e}=\left(x + \frac{k_2}{\sqrt{2}}y\right)I_-+\left(x - \frac{k_1}{\sqrt{2}}y\right)I_+.\]
\end{proof}
\begin{lem0} \label{lem3}
\textit{A function }$f:\,B_c\to A_c$\textit{, }$0<c<1$\textit{, is differentiable if and only if it can be represented as follows}
\[f\left(w\right)=\alpha \left(w_1\right)I_-\,+\,\beta \left(w_2\right)I_+,\]
\textit{where }$w_1= x_1 + {\mathbf i}\,y_1$, $x_1 = x$, $y_1 = - {\mathbf i}\,\frac{k_2}{\sqrt{2}} y$,
$w_2= x_2 + {\mathbf i}\,y_2$, $x_2 = x$, $y_2 = {\mathbf i}\,\frac{k_1}{\sqrt{2}} y$
\textit{and }$\alpha \left(w_1\right)$, $\beta \left(w_2\right)$\textit{ are analytical functions of variables }
$w_1$\textit{, }$w_2,$\textit{ respectively, as follows}
\[\frac{\partial }{\partial y_1}\alpha \left(w_1\right)={\mathbf i} \frac{\partial }{\partial x}\alpha \left(w_1\right),\quad
\frac{\partial }{\partial y_2}\beta \left(w_2\right)={\mathbf i} \frac{\partial }{\partial x}\beta \left(w_2\right).\]
\end{lem0}
\begin{proof}
The sufficiency can be verified directly. Indeed,
\[\frac{\partial }{\partial y}f\left(w\right)=\frac{\partial }{\partial y_1}\alpha \left(w_1\right)
\frac{\partial y_1}{\partial y}I_-+\frac{\partial }{\partial y_2}\beta \left(w_2\right)\frac{\partial y_2}{\partial y}I_+ \]
\[\qquad =-{\mathbf i}\,\frac{k_2}{\sqrt{2}} \frac{\partial }{\partial y_1}\alpha \left(w_1\right)I_-
+ {\mathbf i}\frac{k_1}{\sqrt{2}}\ \frac{\partial }{\partial y_2}\beta \left(w_2\right)I_+\]
\[\qquad = \frac{k_2}{\sqrt{2}} \frac{\partial }{\partial x}\alpha \left(w_1\right)I_- - \frac{k_1}{\sqrt{2}} \frac{\partial }{\partial x}\beta \left(w_2\right)I_+.\]

On the other hand, taking into account Eqs. \eqref{GrindEQ__5_}, \eqref{GrindEQ__6_}, we have
\[{\mathbf e}\frac{\partial }{\partial x}f\left(w\right)=\left(\frac{ k_2}{\sqrt{2}}I_- -
\frac{ k_1}{\sqrt{2}}I_+\right)\left(\frac{\partial }{\partial x}\alpha \left(w_1\right)I_-+ \frac{\partial }{\partial x}\beta \left(w_2\right)I_+\right)\]
\[\quad = \frac{k_2}{\sqrt{2}} \frac{\partial }{\partial x}\alpha \left(w_1\right)I_- - \frac{k_1}{\sqrt{2}} \frac{\partial }{\partial x}\beta \left(w_2\right)I_+.\]

Hence,
\[{\mathbf e}\,\frac{\partial }{\partial x}f\left(w\right)=\frac{\partial }{\partial y}f\left(w\right).\]

Now let us prove necessity. Suppose that a function
\[f\left(w\right)=u_1\left(x,y\right)+ {\mathbf i}\,u_2\left(x,y\right)+{\mathbf e}\,u_3\left(x,y\right)
+{\mathbf i}\,{\mathbf e}\,u_4\left(x,y\right)\]
is monogenic on $B_c$, i.e., ${\mathbf e}\,\frac{\partial }{\partial x}f\left(w\right)=\frac{\partial }{\partial y}f\left(w\right)$.

By using Eq. \eqref{GrindEQ__6_} we can represent $f\left(w\right)$ in the following manner
\[f\left(w\right)=\alpha \left(w_1\right)I_-+\beta \left(w_2\right)I_+,\]
where
\[\alpha \left(w_1\right)=u_1\left(x,y\right)+\frac{ k_2}{\sqrt{2}}u_3\left(x,y\right)
+ {\mathbf i}\,\left( u_2\left(x,y\right) + \frac{ k_2}{\sqrt{2}}u_4\left(x,y\right) \right),\]
\[\beta \left(w_2\right)=u_1\left(x,y\right)-\frac{ k_1}{\sqrt{2}}u_3\left(x,y\right)
+ {\mathbf i}\,\left( u_2\left(x,y\right) - \frac{ k_1}{\sqrt{2}}u_4\left(x,y\right) \right).\]

Consider
\[{\mathbf u}\frac{\partial }{\partial y}f=\frac{\partial }{\partial y_1}\alpha \left(w_1\right)\frac{\partial y_1}{\partial y}I_-
+\frac{\partial }{\partial y_2}\beta \left(w_2\right)\frac{\partial y_2}{\partial y}I_+\]
\[\qquad = -\frac{{\mathbf i} \,k_2}{\sqrt{2}} \frac{\partial \alpha }{\partial y_1}I_-+\frac{{\mathbf i} \,k_1}{\sqrt{2}}\frac{\partial \beta }{\partial y_2}I_+.\]

Then, taking into account Eq. \eqref{GrindEQ__6_}, we have
\[{\mathbf e}\,\frac{\partial }{\partial x}f\left(w\right)=\left(\frac{ k_2}{\sqrt{2}}I_ -
- \frac{ k_1}{\sqrt{2}}I_+\right)\left(\frac{\partial \alpha }{\partial x}I_-+\frac{\partial \beta }{\partial x}I_+\right)\]
\[\quad = \frac{k_2}{\sqrt{2}} \frac{\partial \alpha }{\partial x}I_-
- \frac{k_1}{\sqrt{2}}\frac{\partial \beta }{\partial x} I_+.\]

Therefore,
\[\frac{\partial }{\partial y_1}\alpha \left(w_1\right)= {\mathbf i}\,\frac{\partial }{\partial x}\alpha \left(w_1\right),\]
\[\frac{\partial }{\partial y_2}\beta \left(w_2\right)= {\mathbf i}\,\frac{\partial }{\partial x}\beta \left(w_2\right).\]

Suppose ${\alpha }\left({\omega }_1\right) = {\alpha }_1\left({\omega }_1\right)+{\mathbf i}\,{\alpha }_2\left({\omega }_1\right)$ and ${\beta }\left({\omega }_2\right) = {\beta}_1\left({\omega }_2\right)+{\mathbf i}\,{\beta }_2\left({\omega }_2\right)$.
It follows from the proof of Lemma \ref{lem3} that ${\alpha }_i\left({\omega }_1\right)+{\beta }_j\left({\omega }_2\right)$, $i,j\in \left\{1,2\right\}$  are solutions of Eq. \eqref{GrindEQ__1_} for $0<c<1$.
\end{proof}

\begin{theor10} \label{the2}
$u\left(x,y\right)$\textit{ is a solution of Eq. \eqref{GrindEQ__1_} for }$0<c<1$\textit{ if and only if for some }$i,j\in \left\{1,2\right\}$\textit{ it can be represented as follows}

\[u\left(x,y\right)={\alpha }_i\left({\omega }_1\right)+{\beta }_j\left({\omega }_2\right),\] 

\textit{where }${\alpha }_i\left({\omega }_1\right),{\beta }_j\left({\omega }_2\right)$\textit{ are components of }$\ \alpha \left({\omega }_1\right)$\textit{ and }$\beta \left({\omega }_2\right)$\textit{ of monogenic function }$g\left(\omega \right)$\textit{ in the decomposition \eqref{GrindEQ__15_} i.e.,}

\[f\left(\omega \right)=\alpha \left({\omega }_1\right)I_-+\beta \left({\omega }_2\right)I_+,\] 

\textit{where }$\alpha \left({\omega }_1\right)$\textit{, }$\beta \left({\omega }_2\right)$\textit{ are complex analytical functions of respective variables}.
\end{theor10}

\begin{proof}  As mentioned above  ${\alpha }_i\left({\omega }_1\right)+{\beta }_j\left({\omega }_2\right)$, $i,j\in \left\{1,2\right\}$  are solutions of Eq. \eqref{GrindEQ__1_} for $0<c<1$.

If $u\left(x,y\right)$ is a solution of Eq. \eqref{GrindEQ__1_} for $0<c<1$ much in the same way as in proving Theorem \ref{the1} we can show that $u\left(x,y\right)={\alpha }_i\left({\omega }_1\right)+{\beta }_j\left({\omega }_2\right)$. 
\end{proof}

\end{document}